\documentstyle[12pt]{article}

\def\thebibliography#1{\section*{{\normalsize \begin{center} \bf 
References \end{center} }}\list
  {[\arabic{enumi}]}{\settowidth\labelwidth{[#1]}\leftmargin\labelwidth
    \advance\leftmargin\labelsep
    \usecounter{enumi}}
    \def\newblock{\hskip .11em plus .33em minus -.07em}
    \sloppy
    \sfcode`\.=1000\relax}

\newtheorem{th}{Theorem}

\newtheorem{lem}[th]{Lemma}

\newcommand{\nc}{\newcommand}

\newtheorem{definition}[th]{Definition}
\newtheorem{prop1}[th]{Proposition}
\newtheorem{lemma}[th]{Lemma}
\newtheorem{remark}[th]{Remark}
\newtheorem{cor}[th]{Corollary}
\newtheorem{other}[th]{}
\newtheorem{example}[th]{Example}

\nc{\bother}{\begin{other}}
\nc{\eother}{\end{other}}
\nc{\bex}{\begin{example}}
\nc{\eex}{\end{example}}
\nc{\be}{\begin{equation}}
\nc{\ee}{\end{equation}}
\nc{\bea}{\begin{eqnarray}}
\nc{\eea}{\end{eqnarray}}
\nc{\bth}{\begin{th}}
\nc{\eth}{\end{th}}
\nc{\bdf}{\begin{definition}}
\nc{\edf}{\end{definition}}
\nc{\bpr}{\begin{prop1}}
\nc{\epr}{\end{prop1}}
\nc{\blm}{\begin{lemma}}
\nc{\elm}{\end{lemma}}
\nc{\br }{\begin{remark}}
\nc{\er}{\end{remark}}
\nc{\bcor}{\begin{cor}}
\nc{\ecor}{\end{cor}}

\def\lmn#1{\vadjust{\setbox1=\vtop{\hsize 12mm
\parindent=0pt\baselineskip=9pt
\rightskip=4mm plus 4mm#1}
\hbox{\kern-12mm\smash{\raise .5ex\box1}}}}

\def\R{{\bf R}}
\def\Q{{\bf Q}}

\def\v8{\vskip 8pt}

\def\A{{\cal A}}
\def\s{\sigma}

\def\build#1_#2^#3{\mathrel{\mathop{\kern 0pt#1}\limits_{#2}^{#3}}}
\def\mapright#1{\smash{\mathop{\longrightarrow}\limits^{#1}}}

\def\BZ{\bf Z}
\def\clos#1{<#1>}

\def\mbar{\overline \mu}
\def\e{\emptyset}
\def\l{\ell}
\def\s{\sigma}

\def\A{{\mathcal A}}

\def\G{\Gamma}
\def\ch{\check}

\newcommand{\func}[3]{\mbox{${#1}\colon #2 \rightarrow #3$}}

\def\cp{\coprod_{i=1}^{\l} I}
\def\cps{\coprod_{i=1}^{\l} S^1}
\def\RxI{\R^2\times I}
\def\DxI{D^2\times I}

\date{\today}

\begin{document}

\input epsf.sty

\begin{center}
{\Large The Casson-Walker-Lescop Invariant as a Quantum $3$-manifold Invariant}

\v8\v8
by Nathan Habegger \footnote{Research supported in part by the CNRS.  \\
This and related preprints can be obtained by accessing the WEB at 
{\tt http://www.math.sciences.univ-nantes.fr/preprints/}}

\v8
UMR 6629 du CNRS, Universit\'e de Nantes

D\'epartement de Math\'e\-matiques

Nantes, FRANCE

e-mail: {\tt habegger$\char'100$math.univ-nantes.fr}

\v8
and

\v8
Anna Beliakova 

\v8
Department of Mathematics

University of Berne

Berne, SWITZERLAND

e-mail: {\tt beliak$\char'100$math-stat.unibe.ch}

\v8\v8
(July 24, 1997)
\end{center}

\v8
\begin{abstract} Let $Z(M)$ be the 3-manifold invariant of Le, Murakami
and Ohtsuki.  We give a direct computational proof that the degree 1 part
of $Z(M)$ satisfies $Z_1(M) = \frac{(-1)^{b_1(M)}} {2} \lambda_M$, where 
$b_1(M)$
denotes the first Betti number of $M$ and where $\lambda_M$ denotes the Lescop
generalization of the Casson-Walker invariant of $M$.  Moreover, if $b_1(M)=2$,
we show that $Z(M)$ is determined by $\lambda_M$.

\end{abstract}

\vfill\eject

\section{Introduction and Statement of Results}  A {\em universal} invariant of
3-manifolds, denoted by $Z(M)$, taking values in the graded (completed) algebra
of Feynman diagrams, $A(\emptyset)$, was introduced by T. Le, J. Murakami and
T. Ohtsuki, \cite{LMO}.  The map $M\mapsto Z(M)$, when restricted to integral
homology spheres, was shown by Le in \cite{L1} to be the universal invariant of
finite type (in the sense of T. Ohtsuki, \cite {O}).  In particular, this map
induces an isomorphism from the vector space, which has as basis the set of
diffeomorphism classes of oriented homology spheres (completed with respect to
the Ohtsuki filtration), to $A(\emptyset)$. Thus the invariant $Z(M)$ is a rich
source of information, for homology spheres.

\v8
On the other hand, in \cite{H}, $Z(M)$ was computed for $3$-manifolds $M$ whose
first betti number, $b_1(M)$, is greater than or equal to 3.  More precisely, it
was shown that $Z(M)=1$, if $b_1(M)>3$, and that $Z(M)=\Sigma_n \lambda_M^n
\gamma_n$, if $b_1(M)=3$.  Here $\gamma_n$ denotes a certain nontrivial element
of $A_n(\emptyset)$, and $\lambda_M$ denotes the Lescop invariant of $M$,
\cite{Ls}.   Lescop's invariant is a generalization to all oriented 3-manifolds
of the Casson-Walker invariant of (rational) homology 3-spheres.

More recently, S. Garoufalidis and the first author, \cite {GH}, computed the
universal invariant of a 3-manifold $M$, satisfying $H_1(M,\BZ)=\BZ$, in terms
of the Alexander polynomial of $M$.

\v8
In this paper, we fill in the missing gap for $b_1(M)=2$.  Namely, we have the
following:

\begin{th}\label{betti2} There exist non-zero ${\cal H}_n\in
A_n(\emptyset)$, such that for all closed oriented $M$ satisfying $b_1(M)=2$, 
one has $Z(M)=\Sigma_n\lambda_M^n {\cal H}_n$, where $\lambda_M$ denotes the 
Lescop invariant of $M$.
\end{th}

In order to describe the element ${\cal H}_n$ of $A_n(\emptyset)$, we first
recall that in [LMO], maps
$\iota_n\colon A_{n\l+i}(\cps)\rightarrow A_i(\emptyset)$ were defined for all
$i\ge 0$.    (We set $\iota_n=0$ otherwise.)  We denote by
$p_\l\colon A(\cp)\rightarrow A(\cps)$ the quotient
mapping. We set ${\cal H}_0=1$, and ${\cal H}_n=\iota_n
( p_2(\chi({ {H_{12}^{n} } \over {2^n\;n!} })))\in A_n(\emptyset)$,
where $H_{12}$ is the chinese character of degree 3 which is given by the
'H' diagram below.  
\begin{center}
\mbox{\epsfysize=1.5cm\epsffile{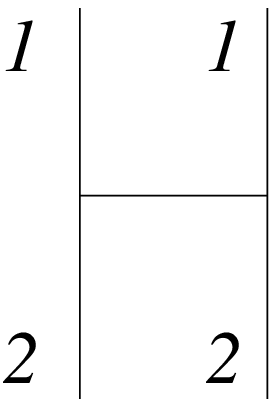}}\\
 Figure 1 {\it  The chinese character $H_{12}$.}
\end{center}


Note that $A_1(\emptyset)$ is 1-dimensional and that
$A_1(\emptyset)^{\otimes n}$ is a direct summand of $A_n(\emptyset)$. 
Moreover, it is easily seen from the definition of $\iota_n$ that the image of 
${\cal H}_n$ in $A_1(\emptyset)^{\otimes n}$ is nonzero.  Hence
${\cal H}_n$ is nonzero.  
\v8

In \cite{HM}, G. Masbaum and the first author investigated the relation
between finite type invariants and Milnor's $\mbar$-invariants of
string-links.   It was shown, for example, that the coefficient of 
$H_{12}$ in (the tree-like part of) the
Kontsevich expansion (the universal finite type invariant of tangles) 
of a (string) link with vanishing linking numbers, is {\em half} the Milnor
invariant $\mu_{1122}$ in the classical index notation.
(In fact, \cite{HM} was motivated by early attempts by the first author to
establish Theorem \ref{betti2}.)

This recent advance makes it possible to prove theorem 1 and to give a direct 
proof of the fact, first
proven in \cite{LMO}, that the linear term in the $LMO$ expansion is, up to a
factor, the Casson-Walker-Lescop invariant.  Namely we prove

\begin{th}\label{lescop} 
$Z_1(M) =  \frac{(-1)^{b_1(M)}} {2} \lambda_M $.
\end{th}



\v8
{\small
\bf Acknowledgements:\rm\ \   
The authors wish to thank Christine Lescop, for her help in establishing
Theorem \ref{lescop}, Dylan Thurston, for suggesting an approach which led to
our proof of Theorem \ref{betti2}, and Gregor Masbaum for his encouragement to
complete this project.  This work began during a visit of the first author to
the University of Berne during June of $1997$.  He thanks Christine Riedtman and
the Berne Mathematics Department for their hospitality.  He also wishes to thank
the CNRS for its support, and the University of Nantes for delegating his time
towards research during the academic year 96-97.  }

\section{The universal finite type invariant}\label{def}

This section provides a brief review of some known facts about the universal
finite type 3-manifold invariant.

\v8
A {\em tangle} in a manifold $M$ is a smooth compact 1-manifold $X$, and
a smooth embedding
\func{T}{(X,\partial X)}{(M^3,\partial M^3)},
transverse to the boundary.
If $X = \cp$, and $M=\DxI$ and the smooth embedding is a fixed
standard  embedding on the boundary, the tangle is called a {\em string link}. 
A {\em framed tangle} is a tangle with a framing which restricts to a
predetermined, fixed framing on the boundary.   (The use of the term framing is
slightly abusive here, as we actually only  require a single non-vanishing
section of the normal bundle.)

Note that tangles may have empty boundary, so that links are special cases of 
tangles. By gluing a string link $\s$,
along the boundary, with the trivial string link, we
have the notion of the {\em closure} $\hat\s$ of
$\sigma$. 

\v8
We refer the reader to~\cite{B1} \cite{B2} for the definition of  the (graded,
completed) $\Q$-coalgebra $\A(X)$ of Feynman diagrams on the 1-manifold $X$. 
Elements of  $\A(X)$ are represented by linear combinations of vertex-oriented
diagrams $X\cup\G$, subject to the AS and IHX  relations.    Here $\G$ is a
uni-trivalent graph, whose univalent vertices lie in the interior of $X$.  It is
customary to refer to the trivalent vertices of $\G$ as {\em internal
vertices} and to the univalent vertices of $\G$ as {\em external
vertices}.
It is also customary to refer to the components of
$X$ as {\em solid},  and to the components of $\G$  as {\em dashed}. 
The space $\A(X)$ is graded by the degree, where the degree of  a diagram is 
half
the number of vertices of $\G$.   The degree $n$ part of $\A(X)$ will be
denoted by $\A_n(X)$.  

We will denote $\A(\cp)$ by $\A(\l)$.  Note that $\A(\l)$ 
has a juxtaposition product (and is thus a Hopf algebra) given by stacking 
diagrams.  An
inclusion of $I$ in $X$ defines an injection of $\A(1)$ into $\A(X)$, and an
action of $\A(1)$ on $\A(X)$.  For
$X=S^1$, such an  inclusion induces an isomorphism of $\A(1)$ with
$\A(S^1)$, which thus inherits a product (and Hopf algebra) structure.
$\A(\e)$ has a product given by disjoint union.

\v8
In~\cite{LM2}, (see also~\cite{LM1}), T. Le and J. Murakami constructed
an  invariant $Z(T) \in \A(X)$, a version of the Kontsevich integral, for
any framed $q$-tangle $\func{T}{X}{\RxI}$.  (This was denoted by $\hat
Z_f(T)$  in~\cite{LM2} and should not be confused with its precursor
$Z_f(T)$.)  

\v8
In \cite{LMO}, Le,
Murakami and Ohtsuki give a 3-manifold invariant, called the {\em LMO
invariant}, defined as follows:
\begin{equation}
Z_n(M) =
\left[\frac{\iota_n(\ch{Z}(L))}{(\iota_n(\ch{Z}(U_{+})))^{\sigma_{+}}
(\iota_n(\ch{Z}(U_{-})))^{\sigma_{-}}}\right]^{(n)} \in
{\A}_n(\e).
\label{eq:LMO}
\end{equation}
Here $\xi^{(n)}$ denotes the degree $n$ part of $\xi \in 
{\A}(\e)$. 
$L\subset S^3$ is a framed link such that surgery on $L$ gives the
3-manifold $M$.  $\ch{Z}(L) = \nu^{\otimes \l}Z(L)$ is obtained by
successively taking the connected sum of $Z(L)$ with $\nu$ along each
component of $L$, where $\nu$ is  the value of $Z$ on the unknot with zero
framing. $U_{\pm}$ is the trivial knot  with $\pm 1$-framing. 
$\sigma_{\pm}$ is the  dimension of the positive and negative eigenspaces of the
linking matrix for the framed link $L\colon \cps \rightarrow {\R}^3$.
 
$$\func{\iota_n}
{{\A}(\cps)}{{\A}(\e)}$$
is a map defined in \cite{LMO}.  The map $\iota_n$, although rather
complicated, is more transparent when evaluated on Chinese characters,
(see below).

\v8
Recall that there is a natural isomorphism $\chi$ of ${\cal A}(\l)$ with ${\cal
B}(\l)$, the 
$\Q$ (Hopf) algebra of so-called {\em Chinese characters}, i.e., uni-trivalent
(dashed) graphs (modulo AS and IHX relations), whose trivalent vertices
are oriented, and whose univalent vertices are labelled by elements of the set
$\{1,\ldots ,\l\}$.
$\chi$ is given by mapping a chinese character to the average of all of the
ways of attaching its labelled edges to the $\l$ strands.  $\chi$ is
comultiplicative, but {\em not} multiplicative.

\v8
One has the following formula, (see \cite {L2}):  For $i \ge 0$, the
composite 
$${\cal B}_{n\l+i}(\l)\mapright {\chi}\A_{n\l+i}(\l) 
\mapright{p_\l}{\A}_{n\l+i}(\cps)\mapright{\iota_n}{\A}_i(\e)$$ 
is zero on those characters which do not have exactly $2n$ univalent vertices of
each label, and on characters which do have exactly $2n$ univalent vertices of
each label is given by the formula
$$\xi \mapsto O_{-2n}(\clos{\xi}),$$
where $\clos{\xi}$ denotes the sum of all ways of joining the univalent vertices
of $\xi$, having the same label, in pairs, and $O_{-2n}$ is the map which sets
circle components equal to $-2n$.

\section{Proof of Theorem 1}

Let $M$ be a 3-manifold obtained by surgery on an algebraically 
split link $L$ in $S^3$ (i.e., $L$ has vanishing pairwise linking numbers).  
We begin by expressing the Lescop invariant of $M$ in terms of the  Milnor
invariants of $L$, assuming that $b_1(M)=2$.  See \cite {C} for geometric
interpretations of the Milnor  invariants $\mu_{ijk}$ and $\mu_{iijj}$ of a
link. 

For  an abelian group $H$, let $|H|$ denote its cardinality if
this is finite and zero otherwise.  

\begin{lemma} 
Let $L$ be an $\l$-component algebraically split link with framing $\mu_{ii}$ on
the $i$-th component, such that $\mu_{11}=\mu_{22}=0$ and $\mu_{ii}\ne 0, i>2.$
Let  $M=S^3(L)$ be the 3-manifold obtained by surgery on $L$.

Then
$$\lambda_M= |\prod^\l_{i=3}\mu_{ii}|\left(\sum^\l_{i=3}\frac{\mu^2_{12i}}
{\mu_{ii}}+\mu_{1122}\right).$$ 
\end{lemma} 
{\bf Proof:}
Recall that in case $b_1(M)=2$ (see \cite{Ls}, T5.2), one has that 
$$\lambda_M= - |{\rm Torsion}(H_1(M))|\;\, {\rm lk}_M(\gamma, \gamma_+).$$
Here, $\gamma$ denotes the framed curve obtained by intersecting two surfaces
whose homology classes generate $H_2(M)$.  $\gamma_+$ denotes the parallel
(using the framing) of $\gamma$. 
Note that $|{\rm Torsion}(H_1(M))|= |\prod^\l_{i=3}\mu_{ii}|$.  Thus it remains to show
that ${\rm lk}_M(\gamma, \gamma_+)=\sum^\l_{i=3}\frac{\mu^2_{12i}}{\mu_{ii}}+\mu_{1122}$.

Since $L$ is algebraically split, we may find Seifert surfaces for each
component avoiding the other components.  The surfaces bounding the first two
components may be capped off in $M$ with disks, to yield closed surfaces
representing the generators of
$H_2(M,\BZ)$. Thus, we must compute ${\rm lk}_M(\gamma, \gamma_+)$, where 
$\gamma$ denotes
the intersection of these two Seifert surfaces.   It follows from
Cochran, \cite{C},
that ${\rm lk}_{S^3}(\gamma, \gamma_+)=-\mu_{1122}$.  

From the geometric interpretation of triple Milnor invariants in terms of the
intersection of three Seifert surfaces, one has that
${\rm lk}_{S^3}(\gamma, L_i)={\rm lk}_{S^3}(\gamma_+, L_i)=-\mu_{12i}$. 
Let $m_i$,
resp. $l_i$, denote the meridian, resp. longitude, of the i-th component. 
Note that since we are doing surgery on the i-th component with framing
$\mu_{ii}$, the cycle  $\mu_{ii}m_i+l_i$ bounds in the complement of
$\gamma_+$ and all other components of $L$.  Therefore one obtains the formula

\begin{eqnarray*}
{\rm lk}_M(\gamma, \gamma_+) & = &
-\sum^\l_{i=3}\mu_{12i}\;{\rm lk}_M(m_i, \gamma_+)+ {\rm lk}_{S^3}(\gamma, 
\gamma_+)\\
& = & 
-\sum^\l_{i=3}\frac{\mu_{12i}}{-\mu_{ii}}\;{\rm lk}_{S^3}(l_i, 
\gamma_+)-\mu_{1122} \\ 
& = & 
-\left(\sum^\l_{i=3}\frac{\mu^2_{12i}}{\mu_{ii}}+\mu_{1122}\right).
\end{eqnarray*} 
$\hfill\Box$
\v8

{\bf Proof of Theorem 1:} 
It is sufficient to compute $Z_n(M)$ for the case when the linking matrix of $L$
is diagonal.  (Indeed, given any symmetric bilinear form over  the
integers, one can find a diagonal matrix $D$ having non-zero determinant, such that the
given form becomes diagonalizable after taking the direct sum with $D$
(see \cite {M}, lemma 2.2). 	  Let $L^\prime$
denote a link whose linking matrix is $D$.
Then $L\sqcup L^\prime$
is equivalent to  an algebraically split link
$L^{\prime\prime}$, via handle sliding (\cite {M}, corollary 2.3).  Let
$M^{\prime\prime}=S^3(L^{\prime\prime})$ and
$M^{\prime}=S^3(L^{\prime})$.  The formula for connected sum (see [LMO])
gives
$$Z_n(M^{\prime\prime})=Z_n(M) |H_1(M^{\prime})|^n.$$
On the other hand, one has that
$$\lambda_{M^{\prime\prime}}=\lambda_M|H_1(M^{\prime})|.$$
Since $|H_1(M^{\prime})|=|\det D|\ne 0$, it follows that
the theorem holds for $M$, provided it holds for $M^{\prime\prime}$.)
\v8

From now on suppose that $L$ is an $\l$-component algebraically split link
with framing $\mu_{ii}$ on the i-th component, such that $\mu_{11}=\mu_{22}=0$
and $\mu_{ii}\ne 0, i>2$.

It follows from \cite{LMO} (which can also be seen from what follows) that the
numerator in the definition of $Z_n(M)$, $\iota_n(\nu^{\otimes \l}Z(L))$,
vanishes in degrees $<n$.  Thus we obtain that 
$$Z_n(M)=(-1)^{n\sigma_+}\;\, \iota_n(\nu^{\otimes \l}Z(L))^{(n)}.$$

One has (see [LM2]) that $Z(L)=p_\l(Z(\sigma)\nu_\l)$, where $\sigma$ is any
string link whose closure is $L$.  Here $\nu_\l\in {\cal A}(\l)$ is obtained from
$\nu=\nu_1$ by the operator  which takes a diagram on the interval to the sum of
all lifts of vertices to  each of the $\l$ intervals. 
 
We claim that

$$
\iota_n(p_\l(\nu^{\otimes \l}Z(L)))^{(n)}
=
\iota_n(p_\l(Z(\sigma)\nu_\l\nu^{\otimes \l}))^{(n)}
$$
$$
=\sum_{m_i \atop m_3+ \ldots +m_\l +m_H =n}
\left( \prod^\l_{k=3} \frac{\mu^{2m_k}_{12k}}{(2m_k)!}
 \frac{(\mu_{kk}/2)^{n-m_k}}{(n-m_k)!}\right) 
 \frac{(\mu_{1122}/2)^{m_H}}{m_H!}
\iota_n(p_\l(\chi(\xi_{m_3,\dots,m_\l,m_H}))),
$$
where $\xi_{m_3,\dots,m_\l,m_H}=W_{123}^{2m_3}\cdots
W_{12\l}^{2m_\l}H_{12}^{m_H}  I_3^{n-m_3}\cdots I_\l^{n-m_\l}$,
and where the chinese caracters 
$I_i$, $W_{ijk}$ and $H_{ij}$ are drawn below.

\begin{center}
\mbox{\epsfysize=1.7cm\epsffile{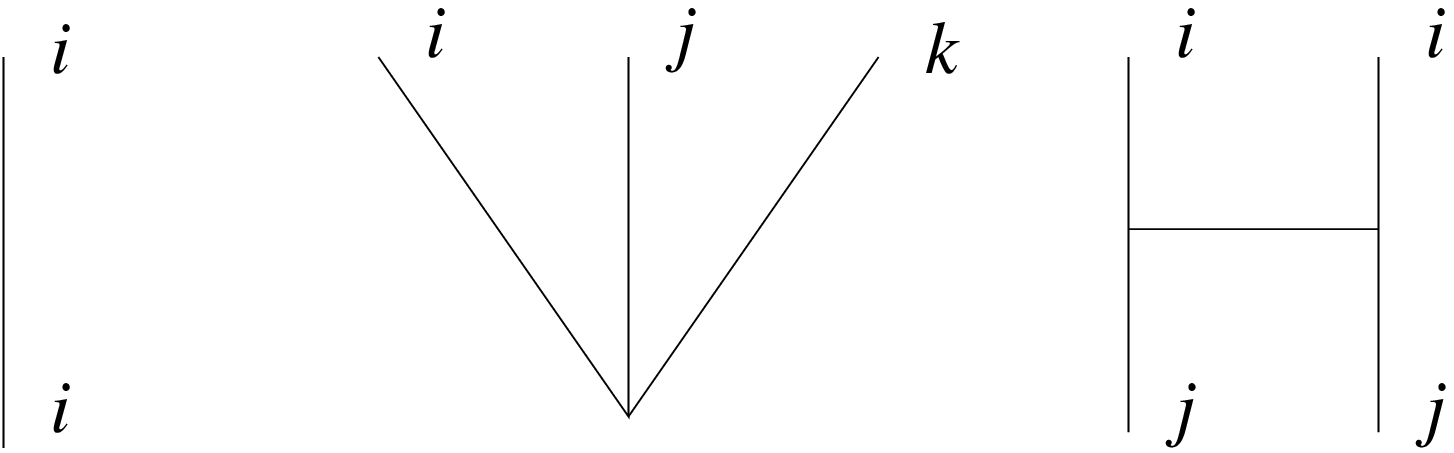}}
\end{center}

To see this, note that 
the coefficients of the chinese caracters 
$I_i$, $W_{ijk}$ and $H_{ij}$ in 
the Kontsevich integral are
$\mu_{ii}/2$, $-\mu_{ijk}$ and $\mu_{iijj}/2$,
respectively, (see \cite{HM}). (N.b., since
$\mu_{ij}=0$, $i\ne j$, the Milnor invariants $\mu_{ijk}$ and 
$\mu_{iijj}$ are well defined integers.)
It follows that one has 
$$ \chi^{-1}(Z(\sigma))= \exp\left( \sum^\l_{k=3}(\mu_{ii}/2)
I_i-\sum^\l_{k=3} \mu_{12k} W_{12k} + (\mu_{1122}/2) H_{12} +
\epsilon\right).$$ (The chinese
character $I_{ij}$, i.e. the interval whose vertices are labelled with $i$
and $j$, does not appear in the above expansion, since $\mu_{ij}=0$, $i\ne
j$.)   Here the exponential is taken
with respect to the disjoint union multiplication in $B(\l)$.  $\epsilon$
consists of terms which are either not trees, are trees of degree $>3$, or are
trees of degree $3$ involving some index other than $1$ and $2$.  These terms
will not be involved in the computation, as they contribute too many vertices
(see below).   A similar formula holds for 
$\chi^{-1}(Z(\sigma)\nu_\l\nu^{\otimes \l})$.

The result claimed above now follows easily, upon remarking
that the chinese characters $\xi_{m_3,\dots,m_\l,m_H}$, with 
$m_3+ \ldots +m_\l +m_H =n$,
are the only terms from the exponential which have $2n$
vertices of each label and have at most $2n$ internal vertices.
(Indeed, since the characters $I_{ij}$, $I_1$ and $I_2$ do not appear, the
only way to obtain $2n$ vertices of label $1$ and $2$, and at most $2n$
internal vertices is to use the characters $W_{12k}$ and $ H_{12}$. Moreover,
using these will yield exactly $2n$ internal vertices.  
Since no other internal vertices are then allowed, the remaining vertices
labelled $k$ must be obtained using the character $I_k$, and since the vertices
of $I_k$ come in pairs, it follows that $W_{12k}$ must occur an even number of
times.)

\v8
Denote by $j^k_n:{\cal B}(k)\to {\cal B}(k-1)$
the map which is zero, if there are not exactly $2n$ univalent vertices labelled
$k$, and otherwise is given by the sum over all
$(2n-1)!!$ ways of joining the $2n$ univalent vertices labelled
$k$ in pairs, and then replacing every circle component by $-2n$. Then
$\iota_n(p_\l(\chi(x)))=j^1_n\circ\cdots\circ j^\l_n
\;(x)$, for $x\in {\cal B}(\l)$.

The theorem is an immediate consequence of the following formula
\be 
\label{3}\frac{j^3_n\circ j^4_n\circ \cdots \circ j^\l_n\;
(\xi_{m_3,\dots,m_\l,m_H})}{
\prod^\l_{k=3}\;2^{n-m_k}\;(n-m_k)!\; (2m_k)!} =\frac{ (-1)^{n\l}}{
\prod^\l_{k=3}\; 2^{m_k}\;m_k!}H^n_{12}. 
\ee
Indeed, using (\ref{3}),  which follows from repeated use of Lemma \ref {WIH}
below,  we obtain that

$$Z_n(M)=(-1)^{n\sigma_+}\;\, \iota_n(\nu^{\otimes \l}Z(L))^{(n)}$$
$$
= (-1)^{n\sigma_+ +n\l}(\prod^\l_{i=3}\mu_{ii})^n
\sum_{m_i \atop m_3+ \ldots +m_\l +m_H =n}
\left( \prod^\l_{k=3} 
\frac{(\mu^2_{12k}/2\mu_{kk})^{m_k}}{(m_k)!}
\right) 
\frac{(\mu_{1122}/2)^{m_H}}{m_H!}
j^1_n\circ j^2_n(H^n_{12})
$$
$$=
\frac{|\prod^\l_{i=3}\mu_{ii}|^n}{2^n \, n!} \left( \sum^\l_{k=3}
\frac{\mu^2_{12k}}{\mu_{kk}} + \mu_{1122}\right)^n 
j^1_n\circ j^2_n(H^n_{12})
$$
$$=
\lambda_M^n
\iota_n ( p_2(\chi({ {H_{12}^{n} } \over {2^n\;n!} }))).
$$
$\hfill\Box$

\begin{lem} 
\label{WIH}
For $i\ge 3$, set $W_{12i}=W$, $H_{12}=H$ and $I_i=I$.  Then
\be\label{12}
\frac{j^i_n\;(W^{2m}I^{n-m})}{(2m)!\;2^{n-m}\; (n-m)!}= \frac{(-1)^n 
}{2^m\;  m!\;}H^m. 
\ee
\end{lem}
\vspace*{0.5cm}

{\bf Proof:} 
Let $J^i_n$ be as in the definition of $j^i_n$, but without the requirement that
we must have exactly $2n$ vertices to be joined.  
To compute $J^i_n(W^{2m}I^k)$, let us fix an $i$-labeled univalent vertex of
$I^k$.  This vertex can either be paired with one of the 
$2m$ $i$-labeled vertices of $W^{2m}$, or with one of the remaining $2k-1$
vertices of $I^k$ (one of which results in a circle component).  Replacing the
circle component by $-2n$, we have that
$$
J^i_n(W^{2m}I^k)=(2m+2k-2-2n)\;J^i_n(W^{2m}I^{k-1}),
$$
and hence inductively we have
\be\label{14}j^i_n(W^{2m}I^{n-m})=(-2)^{n-m}(n-m)!\;J^i_n(W^{2m}).
\ee

Similarly, we have that 
$$
J^i_n(W^{2m})=-(2m-1)\;J^i_n(W^{2m-2})H,
$$
and hence inductively that
\be\label{15}
J^i_n(W^{2m})=(-1)^m(2m-1)!!\; H^m.
\ee

Hence, combining equations (\ref{14}) and (\ref{15}), we obtain 
$$j^i_n(W^{2m}I^{n-m})=(-1)^n 2^{n-m}(n-m)!(2m-1)!!\; H^m,$$
which is equivalent to equation
(\ref{12}), since $(2m-1)!!\;2^m m!=(2m)!$.
$\hfill\Box$

\section{Proof of Theorem 2}

We will denote by $o(n)$ terms of degree greater than or equal to $n$.

In low degrees, the Kontsevich expansion of the
trivial knot can be computed by applying the Alexander-Conway weight system. 
One obtains the formula
$$\chi^{-1}(\nu)= 1+ {\phi_1\over 48} +o(4).$$
From this it follows (see below) that
$$\iota_1(\check{Z}(U_{\pm}))=\mp 1+{\Theta\over 16} + o(2).$$
The graphs $\phi_i$ and $\Theta$ are depicted in Figure 2.
\begin{center}
\mbox{\epsfysize=1.7cm\epsffile{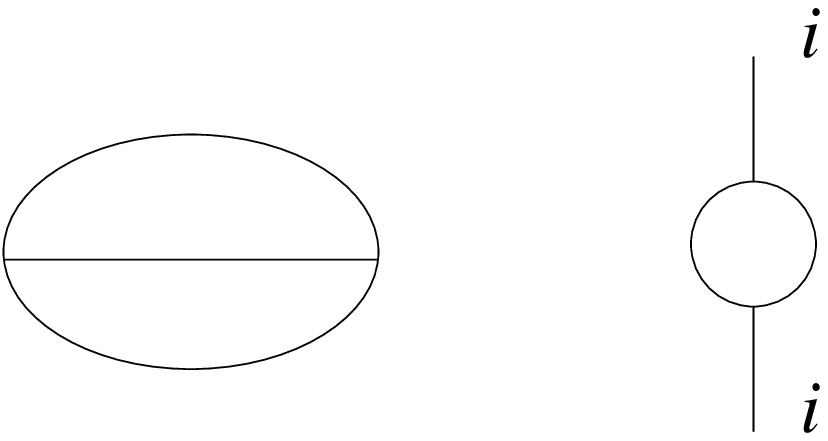}}\\
Figure 2 {\it The graph $\Theta$ and the chinese character $\phi_i$.}
\end{center}

\v8
{\bf Proof of Theorem 2:}
Let $M=S^3(L)$, where $L$ is an $\l$-component algebraically split link which is
the closure of a string link $\s$. (The general case reduces to this one by the
same argument as in the proof of Theorem 1, using the connected sum formulas
for $Z_1$ and $\lambda$, see [LMO], [Ls].)
\v8

We wish to calculate 
$$Z_1(M)=\left[\frac{\iota_1(p_\l(Z(\s)\nu_\l\nu^{\otimes \l}))}
{(-1 + {\Theta\over 16})^{\s_+} (1 + {\Theta\over 16})^{\s_-}}\right]^{(1)}.$$
Denote by $\epsilon$ terms in $\chi^{-1}({\rm log}(Z(\s))$ which either
have more than $2$ internal vertices, or have $2$ internal vertices but an odd
number of external vertices with a given label.  (Such terms do not contribute
to the end expression.)  Then we have
$$
Z(\sigma)= \exp\left(\chi\left( \sum^\l_{i=1}{\mu_{ii}\over 2}
I_i-
\sum_{i<j<k} \mu_{ijk} W_{ijk} +\sum_{i<j}{ \mu_{iijj}\over 2} H_{ij}-
\sum_i{a^{(i)}_1\over 2}\phi_i+\epsilon\right)\right).$$

The coefficient $a^{(i)}_1$ can be calculated
by applying the Alexander-Conway weight system (see e.g., \cite {GH}).  One has 
that $a^{(i)}_1={1\over 2}
\Delta^{\prime\prime}(L_i)(1)$, where $\Delta (K)(t)$ denotes
the Alexander polynomial of $K$ (normalized to be symmetric in 
$t$ and $t^{-1}$
and taking the value $1$ at $t=1$).
In the above expansion of the exponential, one uses the juxtaposition product of
${\cal A}(\l)$.  Note that if we use this product to induce a second product
structure on chinese characters via $\chi$, denoted by $\cdot_\times$, then, for
example, one has that
$$I_k  \cdot_\times I_k=I^2_k + {1\over 6}\phi_k. $$

One deduces that
$$\chi^{-1}(Z(\s)\nu_\l\nu^{\otimes \l})=
 \exp\left( \sum^\l_{i=1}{\mu_{ii}\over 2} I_i- 
\sum_{i<j<k} \mu_{ijk} W_{ijk} +\sum_{i<j}{ \mu_{iijj}\over 2} H_{ij}+
\sum_i{b^{(i)}}\phi_i+\epsilon\right),$$
where $${b^{(i)}}=
2\frac{1}{48}+\frac{1}{6}\frac{1}{2!}\left(\frac{\mu_{ii}}{2}\right)^2 -
{a^{(i)}_1\over 2}=
\frac{2+\mu_{ii}^2-24a^{(i)}_1}{48}.$$
Here the exponential is expanded using the disjoint union product
in ${\cal B}(\l)$.
(Here again, $\epsilon$ denotes terms which either
have more than $2$ internal vertices, or have $2$ internal vertices, but an odd
number of external vertices with a given label.)

It follows that one has 
\begin{equation}
\iota_1(\ch{Z}(L))=(-1)^\l\left(\prod^\l_{i=1}\mu_{ii} +c(L)\Theta\right)+ o(2),
\end{equation}
where the coefficient of $\Theta$, $c(L)$, is given by
$$c(L)= -\sum_i{b^{(i)}}\prod_{j\neq i} \mu_{jj}
+\sum_{i<j}\frac{ \mu_{iijj}}{2}\prod_{k\neq i,j}\mu_{kk}
+\sum_{i<j<k}
\frac{\mu_{ijk}^2}{2}
\prod_{t\neq i,j,k}{\mu_{tt}}
.$$

Hence

\begin{eqnarray*}
Z_1(M) & = &
\left[(-1)^{\s_+}\;\left(1+\frac{(\s_+ -\s_-)}{16}\Theta\right)\;
(-1)^\l\left(\prod^\l_{i=1}\mu_{ii} +c(L)\Theta\right)\right]^{(1)}
\\
& = & 
(-1)^{b_1(M)+\;\s_-}\left(\frac{(\s_+
-\s_-)}{16}\prod^\l_{i=1}\mu_{ii}+c(L)\right)\Theta
\\ 
& = & {(-1)^{b_1(M)}\over 2}\lambda_M\Theta,
\end{eqnarray*}

where
$$\lambda_M={|H_1(M)|\;(\s_+-\s_-)\over 8}\;+
(-1)^{\s_-}2c(L). $$

Note that $\lambda_M$ is indeed Lescop's generalization of the Casson-Walker
invariant.  In the special case of
an algebraically split link in the 3-sphere, Lescop's 
surgery formula (\cite{Ls} 1.4.8) is given as follows:

Let $I$ be a subset of $\{1,\dots,\l\}$. Denote by $L_I$ the link obtained from $L$ by 
forgetting the components whose subscripts do not belong to $I$. Then

$$\lambda_M={|H_1(M)|\;(\s_+-\s_-)\over 8}\;+$$
$$(-1)^{\sigma_-}\left[ \sum_i \left(\zeta(L_{\{i\}}) -
\frac{1}{24}(\mu^2_{ii}+1)\right)\prod_{j\neq i}\mu_{jj} + \sum_{I,\; |I|\neq
0,1}\zeta(L_I)
\prod_{j\neq I}\mu_{jj}\right].$$

Note that in the algebraically split case, the values of the $\zeta$ function 
are given as follows.
$$\zeta(L_I)=\left\{\begin{array}{r@{\quad:\quad}l} a^{(i)}_1 
-\frac{1}{24} & I=\{i\}\\
\mu_{iijj}&I=\{i,j\}\\
\mu_{ijk}^2&I=\{i,j,k\}\\
0&|I|>4\end{array} \right.$$ 
This can easily be established using the results in (\cite {Ls}, chapter 5).
$\hfill\Box$

\v8\v8

\end{document}